\documentclass{math11-en}

\usepackage{tikz-cd} 

\begin{document}

\title{Zariski Topology of (Krasner) Hypermodules}
\maketitle

	\begin{center}
		\begin{alphafootnotes}
			{Behnam Afshar$^{1}$, Reza Ameri$^{1}$ and Hashem Bordbar$^{2}$ }\\[2mm]
			{\small $^{1}$ Department of Mathematics, Statistics and Computer Science, Faculty of Science, University of Tehran, Tehran, Iran.\\
			\small $^{2}$ Center for Information Technologies and Applied Mathematics, University of Nova Gorica, Slovenia.	
				\mail{\url{behnamafshar@ut.ac.ir}}\\
            \mail{\url{rameri@ut.ac.ir}}
        \mail{\url{hashem.bordbar@ung.si}}
    }\\
		
		\end{alphafootnotes}
	\end{center}
%%%%%%%%%%%%%%%%%%%%%%%%%%%%%%%%%%%%%%%%%%%%%%%%%%%%%%%%%%%%%%%%%%%%%%%%%%%%%%%%%%%%%%%%%%%%%%%%%%%%%%%%%%%%%%%%%%	

	\fancyhead[CO]{Zariski Topology of (Krasner) Hypermodules}%%%  write short title of the paper
	\fancyhead[CE]{B. Afshar, R. Ameri and H. Bordbar}%%%  write short names
	
	\hrulefill
	
%%%%%%%%%%%%%%%%%%%%%%%%%%%%%%%%%%%%%%%%%      abstract          %%%%%%%%%%%%%%%%%%%%%%%%%%%%%%%%%%%%%%%%%%%%%%%%%		
\begin{abstract}
Let $R$ be a commutative Krasner hyperring with identity and $M$ an
$R$-hypermodule. We consider the prime spectrum
$\operatorname{Spec}_R(M)$, consisting of all prime subhypermodules of
$M$. In this paper, we introduce a Zariski topology on
$\operatorname{Spec}_R(M)$ and investigate its fundamental topological
properties, including irreducibility, connectedness, compactness, and
separation axioms, providing necessary and sufficient conditions for
each. In particular, we prove that $\operatorname{Spec}_R(M)$ is irreducible
if and only if the nilradical of $M$ is a prime hyperideal, and that
it is connected precisely when the underlying hyperring is not a
nontrivial direct product. We further investigate the functoriality of
the Zariski topology in the category of Krasner hyperrings and
hypermodules.

Moreover, using the fundamental relation $\epsilon^{*}$, the smallest
regular relation on $M$ for which the quotient
$\overline{M}=M/\epsilon^{*}$ becomes an
$\overline{R}=R/\Gamma^{*}$-module, we establish a correspondence
between the Zariski topologies of $M$ and $\overline{M}$. We show that
this correspondence is functorial and thereby relate the
hypermodule-theoretic setting to the classical theory of modules over
rings.
\end{abstract}


%%%%%%%%%%%%%%%%%%%%%%%%%%%%%%%%%%%%%%%%%%%%%%%%%%%%%%%%%%%%%%%%%%%%%%%%%%%%%%%%%%%%%%%%%%%%%%
\section{Introduction}    \label{sec1}

Let $M$ be a unitary module over a commutative ring $R$ with identity.
Then $M$ is called a \emph{multiplication $R$-module} if for every
submodule $N$ of $M$, there exists an ideal $I$ of $R$ such that
$N = IM$. This definition agrees with those of Barnard \cite{barnard}
and El-Bast \cite{Elbast}. Recently, prime submodules and various forms of Zariski topologies on $\operatorname{Spec}(M)$, the set of all prime submodules of $M$, have
been studied in several papers (see, for example,
\cite{amer0, amer00}). Ameri \cite{amer0} investigated prime
submodules of multiplication modules and introduced a notion of
product of submodules, which has since been used to study further
properties of such modules. McCasland \cite{mc1} introduced a topology
on $X = \operatorname{Spec}(M)$, called the Zariski topology, where the
closed sets are defined by
\[
V(N) = \{\, P \in X \mid P \subseteq N \,\},
\]
for submodules $N$ of $M$. However, Durvial \cite{duraivel} and Lu
\cite{lu1} use the term ``Zariski topology'' in a different sense.
In this paper, we focus on the Zariski topology of multiplication
modules over commutative rings with identity, generalizing classical
results on $\operatorname{Spec}(R)$ to $\operatorname{Spec}(M)$, and
study the fundamental properties of this topology.

The theory of hyperstructures was introduced by Marty in 1934 at the
8th Congress of Scandinavian Mathematicians \cite{Marty}. He introduced
hypergroups as a generalization of groups. Since then, this theory has
developed into a rich area of modern algebra. Hyperrings were
introduced by Krasner \cite{Kras1, Kras2}, and later further developed
by many authors. Hypermodules over a hyperring generalize classical
modules over rings. Several types of hyperrings and hypermodules have
been introduced and studied extensively in the literature (see
\cite{amer1, amer8, Bord1, Bord2, bib11, Kras2, mas, Voug2}).

Fixing a Krasner hyperring $R$, one of the most important and well-behaved
classes of $R$-hypermodules is the class of Krasner $R$-hypermodules,
which we simply call \emph{hypermodules} in this paper.

Fundamental relations are among the most important tools in the theory of algebraic hyperstructures, as they allow the passage from hyperstructures to ordinary algebraic structures. For an introduction to the fundamental relations $\beta^{\ast}$ and $\gamma^{\ast}$ on hypergroups and their historical background, see \cite{bib10}. For hyperrings and hypermodules, the fundamental relations $\Gamma^{\ast}$ and $\theta^{\ast}$ are studied in \cite{bib11} and \cite{bib4}, respectively.

In \cite{afshar1}, it is shown that the relations $\beta^{\ast}$ and $\gamma^{\ast}$ are related by
\begin{center}
$\gamma^{\ast} = \delta \ast \beta^{\ast}$,
\end{center}
where $\delta$ is the congruence relation associated with the commutator subgroup. Moreover, \cite{bib2} establishes a one-to-one correspondence between the lattice of strongly regular relations on a regular hypergroup $H$ and the lattice of normal subhypergroups of $H$ containing $\omega_H$. In \cite{bib23}, it is shown that certain sequences of strongly regular relations define a functor, leading to the construction of free objects and tensor products in the category of canonical hypergroups.

In this paper, we investigate the Zariski topology of Krasner hypermodules from the viewpoint of regular and strongly regular relations, which play an important role in the theory of sheaves of
hypermodules. In particular, we make use of the results in \cite{bib25} to study these structures. Our main goal is to apply the topology of relations in the study of sheaves of hypermodules and to extend several known results to a wider class of hypermodules.
%%%%%%%%%%%%%%%%%%%%%%%%%%%%%%%%%%%%%%%%%%%%%%%%%%%%%%%%%%%%%%%%%%%%%%%%%%%%%%%%%%%%%%%%%%%%%%
\section{Preliminaries}

In this section, we recall some basic notions and results on
hyperstructures, including Krasner hyperrings and hypermodules, which
will be used throughout the paper.

Let $H$ be a nonempty set, and let ${P}^*(H)$ denote the family of all
nonempty subsets of $H$. A map
\[
\circ : H \times H \longrightarrow P^*(H)
\]
is called a \emph{hyperoperation}. The hyperoperation $\circ$ can be
naturally extended to nonempty subsets of $H$ by defining, for
$A,B \subseteq H$,
\[
A \circ B = \bigcup_{a \in A,\; b \in B} (a \circ b).
\]

\begin{definition}\label{def1}
A non-empty set $H$ equipped with a commutative hyperoperation $+$ is
called a \emph{canonical hypergroup} if the following conditions hold:
\begin{enumerate}
\item $(x+y)+z = x+(y+z)$ for all $x,y,z \in H$;
\item there exists an element $0 \in H$ such that $x+0 = x$ for all
$x \in H$;
\item for every $x \in H$, there exists a unique element $x' \in H$
such that $0 \in x + x'$ (we denote $x'$ by $-x$ and write
$x-y := x+(-y)$);
\item if $z \in x+y$, then $x \in z-y$ and $y \in z-x$, for all
$x,y,z \in H$.
\end{enumerate}
\end{definition}

\begin{definition}\label{def2}
A non-empty set $R$ equipped with a hyperoperation $+$ and an internal
multiplication $\cdot$ is called a \emph{(Krasner) hyperring} if the
following conditions hold:
\begin{enumerate}
\item $(R,+)$ is a canonical hypergroup;
\item $(R,\cdot)$ is a semigroup;
\item $z(x+y)=zx+zy$ and $(x+y)z=xz+yz$ for all $x,y,z \in R$;
\item $0 \cdot x = 0$ for all $x \in R$.
\end{enumerate}
\end{definition}

A hyperring $R$ is called \emph{commutative} if $xy = yx$ for all
$x,y \in R$. It is called \emph{unitary} if there exists an element
$1 \in R$ such that $1x = x1 = x$ for all $x \in R$.
A commutative unitary hyperring $R$ is called a \emph{hyperfield} if
every nonzero element of $R \setminus \{0\}$ is multiplicatively
invertible.
A nonempty subset $I$ of a hyperring $R$ is called a \emph{left
hyperideal} if $I$ is a subhypergroup of $(R,+)$ and $r x \in I$ for
all $r \in R$ and $x \in I$. Right hyperideals are defined similarly.
A hyperideal is a subset that is both a left and right hyperideal.
A commutative hyperring $R$ is called a \emph{multiplication hyperring}
if, for any hyperideals $A \subsetneq B$ of $R$, there exists a
hyperideal $C$ of $R$ such that $A = BC$.

\begin{definition}\label{def3}
A left hypermodule over a unitary hyperring $R$ is a canonical
hypergroup $M$ equipped with an external operation
$R \times M \to M$, $(r,m) \mapsto rm$, satisfying:
\begin{enumerate}
\item $r(m+n)=rm+rn$ for all $r \in R$ and $m,n \in M$;
\item $(r+s)m = rm + sm$ for all $r,s \in R$ and $m \in M$;
\item $(rs)m = r(sm)$ for all $r,s \in R$ and $m \in M$;
\item $1m = m$ and $0m = 0$ for all $m \in M$.
\end{enumerate}
\end{definition}

An $R$-hypermodule $M$ is a hypermodule over a hyperring $R$.
Throughout this paper, all hypermodules are assumed to be unitary over
a fixed unitary hyperring $R$.
Let $M$ be an $R$-hypermodule. A nonempty subset $N \subseteq M$ is a
\emph{subhypermodule} if it is a hypermodule under the induced
operations. For a hyperideal $I$ of $R$, the set
\[
IM = \bigcup \left\{ \sum_{i=1}^n r_i x_i \,\middle|\,
r_i \in I,\ x_i \in M \right\}
\]
is a subhypermodule of $M$.

Let $R$ and $S$ be hyperrings. A map $f : R \to S$ is called a
\emph{good homomorphism} if, for all $a,b \in R$,
\[
f(a+b)=f(a)+f(b), \quad f(ab)=f(a)f(b), \quad f(0)=0.
\]

An equivalence relation $\sigma$ on a Krasner hyperring $R$ is called
\emph{regular} if
\[
a \,\sigma\, b,\ c \,\sigma\, d \Rightarrow
a+c \,\bar{\sigma}\, b+d,\quad ac \,\sigma\, bd,
\]
and \emph{strongly regular} if
\[
a \,\sigma\, b,\ c \,\sigma\, d \Rightarrow
a+c \,\bar{\bar{\sigma}}\, b+d,\quad ac \,\sigma\, bd.
\]

The fundamental relation $\Gamma^{*}$ on a hyperring $R$ is the smallest
strongly regular relation on $R$. For details, see
\cite{bib11,Voug1}. A relation $\tau$ on an $R$-hypermodule $M$ is
called \emph{regular} if, for all $x,y,z \in M$ and $a \in R$,
\[
x \,\tau\, y \Rightarrow x+z \,\bar{\tau}\, y+z,\quad ax \,\tau\, ay,
\]
and \emph{strongly regular} if
\[
x \,\tau\, y \Rightarrow x+z \,\bar{\bar{\tau}}\, y+z,\quad ax \,\tau\, ay.
\]

The fundamental relation $\theta^{*}$ is the smallest strongly regular
relation on $M$, and the quotient $M/\theta^{*}$ is an
$R/\Gamma^{*}$-module (see \cite{bib4}).

For a commutative hyperring $R$ with identity, a proper hyperideal
$P \subset R$ is called \emph{prime} if $ab \in P$ implies
$a \in P$ or $b \in P$. Every maximal hyperideal is prime.

\begin{proposition}[\cite{bib11}]\label{ii-6}
Let $R$ be a commutative hyperring with identity and $I$ a proper
hyperideal of $R$. Then:
\begin{enumerate}
\item $I$ is prime if and only if $R/I$ is a hyperdomain;
\item $I$ is maximal if and only if $R/I$ is a hyperfield.
\end{enumerate}
\end{proposition}

Recall that a topological space $T$ is:
\begin{itemize}
\item $T_0$ if any two distinct points are topologically distinguishable;
\item $T_1$ if each point is closed;
\item $T_2$ (Hausdorff) if any two distinct points can be separated by
disjoint open sets.
\end{itemize}

Let $\rho$ be a relation on a hypergroup $H$ and $\sigma$ a relation
on $H/\rho$. Then the relation $\sigma * \rho$ on $H$ is defined by
\[
\sigma * \rho = \{(a,b) \in H \times H \mid (\rho(a),\rho(b)) \in \sigma\}.
\]

If $\mathcal{SR}(H)$ denotes the set of strongly regular relations on
a regular hypergroup $H$, and $N(\omega_H)$ the set of normal
subhypergroups containing $\omega_H$, then the map
\[
\varphi : \mathcal{SR}(H) \to N(\omega_H), \quad
\rho \mapsto S_\rho = \{x \in H \mid \rho(x)=e_{H/\rho}\},
\]
is an isomorphism of complete lattices \cite{bib24,bib2}.

%%%%%%%%%%%%%%%%%%%%%%%%%%%%%%%%%%%%%%%%%%%%%%
\section{Zariski topology}\label{sec3}

\begin{definition}\label{v-1}
Let $R$ be a hyperring and let $M$ be a left $R$-hypermodule. A
\emph{prime subhypermodule} of $M$ is a proper subhypermodule
$N \subsetneq M$ such that for any hyperideal $I$ of $R$ and any
subhypermodule $S$ of $M$, the condition $IS \leqslant N$ implies
$S \leqslant N$ or $IM \leqslant N$.
A \emph{primary subhypermodule} of $M$ is a proper subhypermodule
$N \subsetneq M$ such that whenever $r \in R$ and $m \in M$ satisfy
$rm \in N$ and $m \notin N$, there exists $n \in \mathbb{N}$ such that
$r^{n}M \subseteq N$.
\end{definition}

It is clear that every prime subhypermodule is primary. A
\emph{prime hypermodule} $M$ is a hypermodule for which
$\{0_M\}$ is a prime subhypermodule. The set of all prime
subhypermodules of $M$ is called the \emph{prime spectrum} of $M$ and
is denoted by $\operatorname{Spec}(M)$. We also denote by
$\operatorname{mSpec}(M)$ the set of all maximal subhypermodules of
$M$.

\begin{remark}\label{v-2}
For any subhypermodule $N$ of an $R$-hypermodule $M$, define
\begin{center}
$(N:M)=\{\, r \in R \mid rM \subseteq N \,\}$ \ \ \ and \ \ \ 
${Ann}_R(N)=\{\, r \in R \mid rn=0 \text{ for all } n \in N \,\}$.
\end{center}
Both $(N:M)$ and $\operatorname{Ann}_R(N)$ are hyperideals of $R$.
Moreover, $M$ is a prime hypermodule if and only if
$\operatorname{Ann}_R(N)=\operatorname{Ann}_R(M)$ for every nonzero
subhypermodule $N$ of $M$.
Also, a subhypermodule $N$ of $M$ is prime if and only if $M/N$ is a prime hypermodule.
\end{remark}

\begin{example}\label{E1}
\begin{enumerate}
\item[(i)] Let $R=\{0,1\}$ be equipped with the following hyperoperation $+$ and multiplication $\cdot$:
\begin{center}
\begin{tabular}{c|cc}
$+$ & $0$ & $1$\\
\hline
$0$ & $0$ & $1$\\
$1$ & $1$ & $R$
\end{tabular}
\hspace{2cm}
\begin{tabular}{c|cc}
$\cdot$ & $0$ & $1$\\
\hline
$0$ & $0$ & $0$\\
$1$ & $0$ & $1$
\end{tabular}
\end{center}
Then $(R,+,\cdot)$ is a Krasner hyperring whose only hyperideals are
$\{0\}$ and $R$.

Let $M=\{0,a,b\}$. Then $M$ is a Krasner $R$-hypermodule with
hyperoperation $+$ given by:
\begin{center}
\begin{tabular}{c|ccc}
$+$ & $0$ & $a$ & $b$\\
\hline
$0$ & $0$ & $a$ & $b$\\
$a$ & $a$ & $a$ & $M$\\
$b$ & $b$ & $M$ & $b$
\end{tabular}
\end{center}
and external scalar multiplication $R \times M \to M$ defined by
\[
0 \cdot m = 0, \qquad 1 \cdot m = m \quad \text{for all } m \in M.
\]
The subhypermodules of $M$ are $\{0\}$, $\{0,a\}$, $\{0,b\}$, and $M$.
Moreover,
\begin{center}
$\operatorname{Spec}(R)=\{\{0\}\}$, 
${Spec}(M)=\big\{\{0\},\{0,a\},\{0,b\}\big\}$ \  and \ 
${mSpec}(M)=\big\{\{0,a\},\{0,b\}\big\}$.
\end{center}
Therefore, $M$ is a prime Krasner $R$-hypermodule and has no primary
subhypermodule. For every $N \in \operatorname{Spec}(M)$, we have
\[
(N:M)=\{0\}.
\]

\item[(ii)] Consider the same canonical hypergroup $M$ as in (i) with
the following multiplication:
\begin{center}
\begin{tabular}{c|ccc}
$\cdot$ & $0$ & $a$ & $b$\\
\hline
$0$ & $0$ & $0$ & $0$\\
$a$ & $0$ & $a$ & $b$\\
$b$ & $0$ & $a$ & $b$
\end{tabular}
\end{center}
Then $(M,+,\cdot)$ is a Krasner hyperring and hence a Krasner
$M$-hypermodule whose only subhypermodules are $0$ and $M$. Thus,
\[
\operatorname{mSpec}(M)=\operatorname{Spec}(M)=\{0\}.
\]
\end{enumerate}
\end{example}

\begin{proposition}\label{v-3}
Let $R$ be a hyperring and let $M$ be an $R$-hypermodule. Then:
\begin{enumerate}
\item If $N \in \operatorname{Spec}(M)$, then $(N:M) \in \operatorname{Spec}(R)$.
\item If $(N:M) \in \operatorname{mSpec}(R)$, then $N \in \operatorname{Spec}(M)$.
\end{enumerate}
\end{proposition}
\begin{proof}
(1) Let $N \in \operatorname{Spec}(M)$. Since $N \neq M$, it follows that
$(N:M)$ is a proper subhyperring of $R$, hence $(N:M) \neq R$.\\
Let $I$ and $J$ be hyperideals of $R$ such that $IJ \subseteq (N:M)$
and $J \nsubseteq (N:M)$. Then there exist $y \in J$ and $m \in M$
such that $ym \notin N$. Let $x \in I$ and set $z = ym$.\\
Since $IJ \subseteq (N:M)$, we have $xRy \subseteq (N:M)$, and hence
$xRz \subseteq N$. It follows that
\[
(RxR)(Rz) \subseteq N.
\]
Since $N$ is a prime subhypermodule of $M$, we obtain
either $Rz \subseteq N$ or $(RxR)M \subseteq N$.
However, $Rz \nsubseteq N$ because $z \notin N$. Hence
$(RxR)M \subseteq N$, which implies $xM \subseteq N$, and therefore
$x \in (N:M)$. This shows that $(N:M)$ is a prime hyperideal of $R$.

\medskip

(2) Let $(N:M) \in \operatorname{mSpec}(R)$. Consider the quotient
hypermodule $M/N$. For any $m \in M$, the coset $m+N$ satisfies
\begin{align*}
\mathcal{O}_{m+N}
&= \{\, r + (N:M) \mid (r + (N:M))(m+N) = N \,\} \\
&= \{\, r + (N:M) \mid rm + N = N \,\} \\
&= \{(N:M)\}.
\end{align*}
Thus $M/N$ is a torsion-free $R/(N:M)$-hypermodule.
Let $I$ be a hyperideal of $R$ and $S$ a subhypermodule of $M$ such
that $S \nsubseteq N$ and $IS \leqslant N$. Let $i \in I$ and
$x \in S \setminus N$. Then
\[
(i + (N:M))(x+N) = ix + N = N.
\]
Since $M/N$ is torsion-free over $R/(N:M)$, it follows that
$i \in (N:M)$, hence $IM \subseteq N$.\\
Therefore $N$ is a prime subhypermodule of $M$.
\end{proof}

\begin{example}
\begin{enumerate}
\item[(i)] Consider the hyperring $M$ in Example \ref{E1} (ii), and define
\begin{center}
$\operatorname{Mat}(M)=\left\{
\begin{bmatrix}
a & b\\
c & d
\end{bmatrix}
\ \middle|\ a,b,c,d \in M
\right\}.$
\end{center}
Then $\operatorname{Mat}(M)$ is a Krasner $M$-hypermodule with the usual
external scalar multiplication and hyperoperation given by
\[
\begin{bmatrix}
a & b\\
c & d
\end{bmatrix}
+
\begin{bmatrix}
a' & b'\\
c' & d'
\end{bmatrix}
=
\left\{
\begin{bmatrix}
x & y\\
z & w
\end{bmatrix}
\ \middle|\ 
x \in a+a',\ y \in b+b',\ z \in c+c',\ w \in d+d'
\right\}.
\]

Let
\[
N=\left\{
\begin{bmatrix}
a & b\\
0 & c
\end{bmatrix}
\ \middle|\ a,b,c \in M
\right\}.
\]
Then
\begin{align*}
(N:\operatorname{Mat}(M))
&= \{ m \in M \mid m \cdot \operatorname{Mat}(M) \subseteq N \}\\
&= \{ m \in M \mid m \cdot x = 0 \ \text{for all } x \in \operatorname{Mat}(M) \}\\
&= \{0\}.
\end{align*}

Since $\operatorname{mSpec}(M)=\{0\}$, it follows that
$N \in \operatorname{Spec}(\operatorname{Mat}(M))$.

\item[(ii)] Let $x \oplus y = \{x,y,x+y\}$ and $x \odot y = \{xy\}$
for all $x,y \in \mathbb{Z}$. Then $(\mathbb{Z},\oplus,\odot)$ is a
Krasner hyperdomain, whose hyperideals are of the form $n\mathbb{Z}$
for $n \in \mathbb{Z}$.

Moreover,
\[
\operatorname{Spec}(\mathbb{Z})
=
\{p\mathbb{Z} \mid p=0 \ \text{or}\ p \ \text{is prime}\},
\quad
\operatorname{mSpec}(\mathbb{Z})
=
\{p\mathbb{Z} \mid p \ \text{is prime}\}.
\]

On the other hand, $\mathbb{Z}$ is a Krasner $\mathbb{Z}$-hypermodule,
and its subhypermodules are $n\mathbb{Z}$ for $n \in \mathbb{Z}$.
For every prime number $p$, we have
\[
(p\mathbb{Z}:\mathbb{Z}) = p\mathbb{Z} \in \operatorname{mSpec}(\mathbb{Z}).
\]
It is also easy to see that $\mathbb{Z}$ is a prime hypermodule.
\end{enumerate}
\end{example}

\begin{theorem}\label{v-4}
Let $M$ be an $R$-hypermodule. Assume that for every $r \in R$
(resp. for every $m \in M$) we have $r-r=\{0\}$ (resp. $m-m=\{0\}$).
Then
\[
\operatorname{mSpec}(M) \subseteq \operatorname{Spec}(M).
\]
\end{theorem}

\begin{proof}
Let $H$ be a canonical hypergroup. It is well known that $H$ is a group
if and only if $h-h=\{0\}$ for every $h \in H$.\\
Let $N \in \operatorname{mSpec}(M)$. Then $M/N$ is a simple
unitary $R$-hypermodule. In particular, for any nonzero element
$m+N \in M/N$, we have
\[
M/N = R(m+N).
\]

Moreover,
\[
(N:M) = \operatorname{Ann}_R(M/N) = \operatorname{Ann}_R(m+N),
\]
which is a hyperideal of $R$.

Let $(N:M) \subsetneq J \subseteq R$ be a hyperideal and let
$j \in J \setminus (N:M)$. Then $jm+N \neq N$, hence
\[
\bar{x} := jm+N \neq N \in M/N.
\]
Since $M/N$ is simple, we have $R\bar{x} = M/N$. Thus there exists
$r \in R$ such that
\[
r\bar{x} = rjm+N = m+N.
\]

Hence
\[
(rjm - m) + N = N,
\]
so
\[
(rj - 1)(m+N) = N,
\]
which implies $1-rj \in \operatorname{Ann}_R(m+N)$ and therefore
$1-rj \in J$.

\medskip

If $m-m=\{0\}$ for all $m \in M$, the same argument gives
\[
(rj-1)(m+N)=N,
\]
and hence $rj-1 \in \operatorname{Ann}_R(m+N) \subseteq J$.

In both cases, since $rj \in J$ and $1-rj \in J$, we conclude that
$1 \in J$, and hence $J = R$.

Therefore $(N:M)$ is a maximal hyperideal of $R$, i.e.
\[
(N:M) \in \operatorname{mSpec}(R).
\]
By Proposition~\ref{v-3}, it follows that
$N \in \operatorname{Spec}(M)$.
\end{proof}

Let $M$ be an $R$-hypermodule and define
\[
T(M)=\{\, V(S)\mid S \leqslant M \,\},
\]
where
\[
V(S)=\{\, N \in \operatorname{Spec}(M) \mid S \subseteq N \,\}.
\]

Then we have:
\[
V(M)=\emptyset, \qquad V(\{0\})=\operatorname{Spec}(M),
\]
and for any family $\{S_i\}_{i \in I}$ of subhypermodules of $M$,
\[
\bigcap_{i \in I} V(S_i) = V\left(\sum_{i \in I} S_i\right).
\]

However, in general, $T(M)$ is not closed under finite unions.

\begin{definition}\label{v-5}
If $T(M)$ is closed under finite unions, then $M$ is called a
\emph{top hypermodule}. In this case, $T(M)$ defines a topology on
$\operatorname{Spec}(M)$, called the \emph{Zariski topology of $M$}.
\end{definition}

\begin{remark}\label{v-6}
Let $R$ be a commutative hyperring with identity. Then the prime
subhypermodules of $R$ coincide with the prime hyperideals of $R$.
Hence $R$ is a top $R$-hypermodule, and $T(R)$ coincides with the
classical Zariski topology on $\operatorname{Spec}(R)$.

Let $M$ be an $R$-hypermodule, $N$ a subhypermodule of $M$, and
$P \in V(N) \cup V(IM)$. Then $IN \subseteq P$, and hence
$P \in V(IN)$. Conversely, if $P \in V(IN)$, then either
$N \subseteq P$ or $IM \subseteq P$. Therefore,
\[
V(N) \cup V(IM) = V(IN).
\]
Similarly, for hyperideals $I$ and $J$ of $R$, we have
\[
V(IM) \cup V(JM) = V(IJM).
\]

If for every subhypermodule $N$ of $M$ there exists a hyperideal $I$
of $R$ such that $N = IM$, then $M$ is called a
\emph{multiplication hypermodule}. Every multiplication hypermodule is
a top hypermodule.
\end{remark}

\begin{example}\label{v-7}
Let $M = V_{4}$ and consider the Krasner hyperring
$(\overline{\mathbb{Z}}, \oplus, \odot)$, where
\[
\overline{\mathbb{Z}} = \{\bar{n} \mid n \in \mathbb{Z}\}, \quad
\bar{n} = \{-n, n\},
\]
with operations defined by
\[
\bar{n} \oplus \bar{m} = \{\overline{n+m}, \overline{n-m}\}, \qquad
\bar{n} \odot \bar{m} = \overline{nm}.
\]

Then $V_{4}$ is a $\overline{\mathbb{Z}}$-Krasner hypermodule. By
Theorem~\ref{v-4}, we have
\[
\operatorname{mSpec}(V_{4}) \subseteq \operatorname{Spec}(V_{4}).
\]

Hence,
\[
\operatorname{Spec}(V_{4}) = \{A,B,C\},
\]
where
\[
A = \{e,a\}, \quad B = \{e,b\}, \quad C = \{e,c\}.
\]

However,
\[
A \wedge B = \{e\}, \qquad
V(A) \cup V(B) = \{A,B\} \neq V(\{e\}).
\]

Therefore, $V(A) \cup V(B) \notin T(V_{4})$, and hence $T(V_{4})$ is not
closed under finite unions.
\end{example}

In general, if $M$ is a hypermodule such that
$\operatorname{mSpec}(M) \subseteq \operatorname{Spec}(M)$ and
$M_{5}$ can be embedded into the lattice of subhypermodules of $M$
with $1_{M_{5}} = M$, then $T(M)$ is not closed under finite unions.

\begin{definition}\label{v-8}
For a subhypermodule $S \leqslant M$, define
\[
W(S) = \operatorname{Spec}(M) \setminus V(S)
= \{\, N \in \operatorname{Spec}(M) \mid S \nsubseteq N \,\}.
\]

Let
\[
W(M) = \{\, W(S) \mid S \leqslant M \,\}.
\]

Since
\[
\operatorname{Spec}(M) = \bigcup_{S \leqslant M} W(S),
\]
the collection $W(M)$ forms a subbasis for a topology on
$\operatorname{Spec}(M)$.

The topology generated by $W(M)$ (i.e., all unions of finite
intersections of elements of $W(M)$) is called the
\emph{Zariski topology on $M$}, denoted by $\tau(M)$.
\end{definition}
%%%%%%%%%%%%%%%%%%%%%%%%%%%%%%%%%%%%%%%%%%%%%%
Here we introduce a special class of top hypermodules.

Recall that an $R$-hypermodule $M$ is called a
\emph{multiplication hypermodule} if for every subhypermodule $N$ of
$M$, there exists a hyperideal $I$ of $R$ such that $N = IM$.
R. Ameri and M. Norouzi in \cite{amer2} introduced and studied
$(m,n)$-multiplication hypermodules, which generalize multiplication
hypermodules (in particular for $m=n=2$). The following result
immediately follows from \cite{amer2} in the case $m=n=2$.

Let $M$ be an $R$-hypermodule and $N$ a subhypermodule of $M$. Then
\[
I_N = \{\, r \in R \mid rM \subseteq N \,\}
\]
is a hyperideal of $R$. Some authors also denote $I_N$ by $[N:M]$.

\begin{theorem}\label{3}
Let $M$ be a multiplication $R$-hypermodule and let $N$ be a
subhypermodule of $M$. Then
\[
N = I_N M.
\]
\end{theorem}

The class of all multiplication $R$-hypermodules is denoted by
$\mathbf{HMod}_m(R)$. Let $\mathbf{I}(R)$ denote the set of all
hyperideals of $R$. Then $(\mathbf{I}(R), \subseteq)$ is a partially
ordered set.

For an $R$-hypermodule $M$, let $\mathbf{Subh}_R(M)$ denote the set
of all subhypermodules of $M$. Then $(\mathbf{Subh}_R(M), \subseteq)$
is also a partially ordered set.

Define the map
\[
\eta_M : \mathbf{I}(R) \longrightarrow \mathbf{Subh}_R(M), \qquad
I \longmapsto IM.
\]
This map preserves inclusion.

Clearly, an $R$-hypermodule $M$ is a multiplication hypermodule if
and only if $\eta_M$ is surjective.

Finally, every multiplication hyperring induces a multiplication
structure on its associated hypermodules, and conversely every
multiplication $R$-hypermodule arises naturally from this setting
under trivial hyperoperations.

\begin{example}\label{3.6}
Let $(P,+,\cdot)$ be a ring with identity, and let $G$ be a normal
subgroup of the multiplicative semigroup $(P,\cdot)$, i.e.,
$xG = Gx$ for all $x \in P$.

Then the quotient set $\bar{P} = P/G$ becomes a hyperring if we define
addition and multiplication by
\[
xG \oplus yG = \{(xp + yq)G \mid p,q \in G\}, \qquad
xG \odot yG = xyG.
\]
It was shown in \cite{mas} that $(\bar{P}, \oplus, \odot)$ is a
hyperring with identity, called the \emph{quotient hyperring} of $P$
by $G$.
It is easy to verify that there is a one-to-one correspondence between
the ideals of $P$ and the hyperideals of $\bar{P}$ given by
\[
I \triangleleft P \quad \longleftrightarrow \quad \bar{I} = I/G \triangleleft \bar{P}.
\]

Now let $\bar{A}$ and $\bar{B}$ be hyperideals of $\bar{P}$ such that
$\bar{A} \subseteq \bar{B}$. Then $\bar{A} = A/G$ and
$\bar{B} = B/G$, where $A,B \triangleleft P$ with $A \subseteq B$.
Since $P$ is a multiplication ring, there exists an ideal $C$ of $P$
such that $A = BC$. Hence
\[
\bar{A} = \bar{B}\,\bar{C},
\]
and therefore $\bar{P}$ is a multiplication hyperring.

\medskip

As a special case, let $F$ be a field and $F^*$ its multiplicative
group. Then
\[
\bar{P} \cong \{0,1\},
\]
where $\mathbf{K} = \{0,1\}$ becomes a hyperfield under usual
multiplication and hyperaddition given by:
\begin{center}
\begin{tabular}{c|cc}
$+$ & $0$ & $1$\\
\hline
$0$ & $0$ & $1$\\
$1$ & $1$ & $\mathbf{K}$
\end{tabular}
\end{center}
In particular, $1+1 = \{0,1\}$, and $\mathbf{K}$ is the classical
\emph{Krasner hyperfield}.
\end{example}

At the following, we present a method to construct a multiplication
hypermodule from a given multiplication module.

\begin{example}\label{3.7}
In Example~\ref{3.6}, assume that the ring $P$ has identity $1 \in P$
with $1 \in G$. We show that if $M$ is a multiplication $P$-module,
then the induced structure $\bar{M}$ is a multiplication
$\bar{P}$-hypermodule.

Let $M$ be a multiplication $P$-module and let $\bar{I}$ be an ideal
of $\bar{P}$. Define
\[
I = \{\, x \in P \mid xG \in \bar{I} \,\}.
\]
Since $1 \in G$, for any $x,y \in I$ we have
\[
(x-y)G \in xG \ominus yG \subseteq \bar{I},
\]
hence $x-y \in I$. Moreover, for any $x \in I$ and $r \in P$, we have
$xG \in \bar{I}$, and thus
\[
rG \odot xG = rxG \in \bar{I},
\]
which implies $rx \in I$. Therefore, $I$ is an ideal of $P$.
Since $M$ is a multiplication $P$-module, there exists an ideal $J$
of $P$ such that $I = JM$.
We now show that $\bar{I} = \bar{J}\,\bar{M}$.

\medskip

\noindent\textbf{Step 1:} $\bar{I} \subseteq \bar{J}\bar{M}$.

Let $xG \in \bar{I}$. Then $x \in I = JM$, so
\[
x = \sum_{i=1}^n r_i m_i,
\]
where $r_i \in J$ and $m_i \in M$. Hence
\[
xG = \left(\sum_{i=1}^n r_i m_i\right)G
\in \bigoplus_{i=1}^n (r_i m_i)G
= \bigoplus_{i=1}^n (r_i G)(m_i G)
\subseteq \bar{J}\bar{M}.
\]
Thus $\bar{I} \subseteq \bar{J}\bar{M}$.

\medskip

\noindent\textbf{Step 2:} $\bar{J}\bar{M} \subseteq \bar{I}$.

Let $xG \in \bar{J}\bar{M}$. Then
\[
xG \in \bigoplus_{i=1}^n (r_i G)(m_i G),
\]
where $r_i G \in \bar{J}$ and $m_i G \in \bar{M}$.
Hence there exist $q_1,\dots,q_n \in G$ such that
\[
x = \sum_{i=1}^n (r_i m_i) q_i
= \sum_{i=1}^n (r_i q_i)m_i \in JM = I.
\]
Therefore $xG \in \bar{I}$, and so $\bar{J}\bar{M} \subseteq \bar{I}$.

\medskip

Consequently,
\[
\bar{I} = \bar{J}\bar{M},
\]
and thus $\bar{M}$ is a multiplication $\bar{P}$-hypermodule.
\end{example}

\begin{theorem}
Every multiplication $R$-hypermodule is a top hypermodule.
\end{theorem}

\begin{proof}
The result follows immediately from the identity
$V(IM)\cup V(JM)=V(IJM)$ for hyperideals $I,J$ of $R$.
\end{proof}

%%%%%%%%%%%%%%%%%%%%%%%%%%%%%%%%%%%%%%%%%%%%%
The \textit{prime dimension} of an $R$-hypermodule $M$, denoted by
$\dim(M)$, is defined as the maximal length of a chain of prime
subhypermodules in $\operatorname{Spec}(M)$.
If every prime subhypermodule of $M$ is maximal in
$\operatorname{Spec}(M)$, then $\dim(M)=0$. If
$\operatorname{Spec}(M)=\emptyset$, we set $\dim(M)=-1$.

\begin{lemma}\label{v-10}
$\operatorname{Spec}(M)$ is a $T_1$-space if and only if
$\dim(M)\leq 0$.
In particular, if $M$ is a prime hypermodule, then
$\operatorname{Spec}(M)$ is a $T_1$-space if and only if
$\operatorname{Spec}(M)=\{(0)\}$.
\end{lemma}

\begin{proof}
$(\Rightarrow)$ Assume that $\operatorname{Spec}(M)$ is a $T_1$-space.
Then every singleton $\{P\}$, where $P \in \operatorname{Spec}(M)$,
is closed.
Hence, for each $P \in \operatorname{Spec}(M)$, there exists a family
$\{V(N_i)\}_{i \in I}$ such that
\[
\{P\} = \bigcap_{i \in I} V(N_i),
\]
where each $N_i \leq M$.
Fix $i \in I$. Since $P \in V(N_i)$, we have $N_i \subseteq P$.
If there existed $Q \in \operatorname{Spec}(M)$ such that
$P \subsetneq Q$, then $N_i \subseteq Q$, hence $Q \in V(N_i)$ for all
$i \in I$. This implies $Q \in \{P\}$, a contradiction.\\
Thus, no prime subhypermodule properly contains $P$, so every
$P \in \operatorname{Spec}(M)$ is maximal. Hence
$\dim(M) \leq 0$.

\medskip

$(\Leftarrow)$ Assume that $\dim(M)\leq 0$. Then every prime
subhypermodule of $M$ is maximal in $\operatorname{Spec}(M)$.
Let $P \in \operatorname{Spec}(M)$. Since $P$ is maximal, we have
\[
V(P)=\{P\}.
\]
Thus every singleton in $\operatorname{Spec}(M)$ is closed, and hence
$\operatorname{Spec}(M)$ is a $T_1$-space.
\end{proof}

\begin{theorem}\label{v-11}
Let $M$ be an $R$-hypermodule and let $N \leq M$. Then
$\operatorname{Spec}(M)$ is a cofinite topological space if and only if
$\dim(M)\leq 0$ and for every subhypermodule $N \leq M$, either
$V(N)$ is finite or $V(N)=\operatorname{Spec}(M)$.
\end{theorem}

\begin{proof}
$(\Rightarrow)$ Assume that $\operatorname{Spec}(M)$ is a cofinite
topological space. Then every closed set is either finite or the whole
space.
Since every cofinite space is $T_1$, by Lemma~\ref{v-10} we obtain
$\dim(M)\leq 0$.\\
Now let $N \leq M$. Since $V(N)$ is closed in $\operatorname{Spec}(M)$,
it follows that $V(N)$ is either finite or equal to
$\operatorname{Spec}(M)$.

\medskip

$(\Leftarrow)$ Conversely, assume that $\dim(M)\leq 0$ and that for
every subhypermodule $N \leq M$, the set $V(N)$ is either finite or
equal to $\operatorname{Spec}(M)$.
Let $F \subseteq \operatorname{Spec}(M)$ be a closed set. Then
$F = V(N)$ for some subhypermodule $N \leq M$. By assumption, $F$ is
either finite or the whole space.
Hence all closed sets in $\operatorname{Spec}(M)$ are finite or
$\operatorname{Spec}(M)$ itself, and therefore
$\operatorname{Spec}(M)$ is a cofinite topological space.
\end{proof}

\begin{theorem}\label{v-111}
Let $R$ be a hyperring and $M$ an $R$-hypermodule. Then
\[
\operatorname{Spec}\!\left(\frac{M}{\theta^{*}}\right)
=
\left\{
\frac{P}{\theta^{*}} \; \middle|\; P \in \operatorname{Spec}(M)
\right\}.
\]
\end{theorem}

\begin{proof}
Let $N$ be a prime subhypermodule of $M/\theta^{*}$. Define
\[
P = \{\, m \in M \mid \theta^{*}(m) \in N \,\}.
\]
Then $N = P/\theta^{*}$.
We first show that $P$ is a subhypermodule of $M$.\\
Let $m_1, m_2 \in P$ and $r \in R$. Since $N$ is a subhypermodule of
$M/\theta^{*}$, we have
\[
\theta^{*}(m_1) - \theta^{*}(m_2) \in N,
\qquad
\Gamma^{*}(r)\,\theta^{*}(m_1) \in N.
\]
Hence
\[
\theta^{*}(m_1 - m_2) \in N,
\qquad
\theta^{*}(rm_1) \in N,
\]
which implies $m_1 - m_2 \in P$ and $rm_1 \in P$. Thus $P$ is a
subhypermodule of $M$.

Since $N \neq M/\theta^{*}$, it follows that $P \neq M$.
We now prove that $P$ is prime. Let $I$ be a hyperideal of $R$ and
$S \leq M$ such that $IS \subseteq P$. Passing to quotients gives
\[
\frac{I}{\Gamma^{*}} \cdot \frac{S}{\theta^{*}}
\subseteq \frac{P}{\theta^{*}} = N.
\]
Since $N$ is prime in $M/\theta^{*}$, we obtain
\[
\frac{S}{\theta^{*}} \subseteq \frac{P}{\theta^{*}}
\quad \text{or} \quad
\frac{IM}{\theta^{*}} \subseteq \frac{P}{\theta^{*}}.
\]
Hence $S \subseteq P$ or $IM \subseteq P$. Therefore $P \in \operatorname{Spec}(M)$.
This shows that every prime subhypermodule of $M/\theta^{*}$ arises
from a prime subhypermodule of $M$.

Conversely, if $P \in \operatorname{Spec}(M)$, then the same argument
shows that $P/\theta^{*}$ is a prime subhypermodule of $M/\theta^{*}$.
Therefore,
\[
\operatorname{Spec}\!\left(\frac{M}{\theta^{*}}\right)
=
\left\{
\frac{P}{\theta^{*}} \; \middle|\; P \in \operatorname{Spec}(M)
\right\}.
\]
\end{proof}

\begin{theorem}\label{v-12}
Let $f:M \rightarrow N$ be a strong homomorphism of $R$-hypermodules.
Then the map
\[
\bar{f}: \operatorname{Spec}(N) \longrightarrow \operatorname{Spec}(M),
\qquad
\bar{f}(P)=f^{-1}(P),
\]
is well-defined and continuous.
\end{theorem}

\begin{proof}
Let $P \in \operatorname{Spec}(N)$. We first show that $f^{-1}(P)$ is a
prime subhypermodule of $M$. Clearly, $f^{-1}(P)$ is a subhypermodule of $M$. Let $I$ be a hyperideal of $R$ and $T \leq M$ such that
\[
IT \subseteq f^{-1}(P).
\]
Applying $f$, we obtain
\[
f(IT) \subseteq P.
\]
Since $f$ is a strong homomorphism, it follows that
\[
If(T) \subseteq P.
\]
Because $P$ is a prime subhypermodule of $N$, we have
\[
f(T) \subseteq P \quad \text{or} \quad IN \subseteq P.
\]
If $f(T) \subseteq P$, then $T \subseteq f^{-1}(P)$.
If $IN \subseteq P$, then in particular $IM \subseteq f^{-1}(P)$.
Hence $f^{-1}(P)$ is a prime subhypermodule of $M$, and therefore
$\bar{f}$ is well-defined.

\medskip

Now let $V(T)$ be a closed subset of $\operatorname{Spec}(M)$, where
$T \leq M$. Then
\begin{align*}
\bar{f}^{-1}(V(T))
&= \{P \in \operatorname{Spec}(N) \mid f^{-1}(P) \in V(T)\} \\
&= \{P \in \operatorname{Spec}(N) \mid T \subseteq f^{-1}(P)\} \\
&= \{P \in \operatorname{Spec}(N) \mid f(T) \subseteq P\} \\
&= V(f(T)).
\end{align*}
Since $V(f(T))$ is closed in $\operatorname{Spec}(N)$, the map
$\bar{f}$ is continuous.
\end{proof}

\begin{proposition}
The assignment
\[
F : Z.\mathrm{Top}.R\text{-HM} \longrightarrow Z.\mathrm{Top}.\frac{R}{\Gamma^{\ast}}\text{-M}
\]
defined by
\begin{itemize}
\item $F(\operatorname{Spec}(M)) = \operatorname{Spec}\!\left(\frac{M}{\theta^{\ast}}\right)$,
\item for every morphism $f : \operatorname{Spec}(M) \to \operatorname{Spec}(N)$,
the induced map
\[
f^{\ast} : \operatorname{Spec}\!\left(\frac{M}{\theta^{\ast}}\right)
\longrightarrow
\operatorname{Spec}\!\left(\frac{N}{\theta^{\ast}}\right)
\]
is given by
\[
f^{\ast}(\theta^{\ast}(x)) = \theta^{\ast}(f(x)),
\]
\end{itemize}
defines a functor.
\end{proposition}
\begin{proof}
We first show that $f^{\ast}$ is well-defined.
If $\theta^{\ast}(x)=\theta^{\ast}(y)$, then $x \,\theta^{\ast}\, y$.
Since $f$ is compatible with the fundamental relation, it follows that
$f(x)\,\theta^{\ast}\, f(y)$, and hence
\[
\theta^{\ast}(f(x)) = \theta^{\ast}(f(y)).
\]
Thus $f^{\ast}$ is well-defined.

\medskip

We now show that $f^{\ast}$ is continuous. Let
\[
V\!\left(\frac{S}{\theta^{\ast}}\right)
=
\left\{
\frac{P}{\theta^{\ast}} \in \operatorname{Spec}\!\left(\frac{N}{\theta^{\ast}}\right)
\;\middle|\;
P \in V(S)
\right\}
\]
be a closed set in $\operatorname{Spec}(N/\theta^{\ast})$.
Then
\begin{align*}
(f^{\ast})^{-1}
\left(
V\!\left(\frac{S}{\theta^{\ast}}\right)
\right)
&=
\left\{
\frac{Q}{\theta^{\ast}} \in \operatorname{Spec}\!\left(\frac{M}{\theta^{\ast}}\right)
\;\middle|\;
f^{\ast}\!\left(\frac{Q}{\theta^{\ast}}\right)
\in V\!\left(\frac{S}{\theta^{\ast}}\right)
\right\} \\
&=
\left\{
\frac{Q}{\theta^{\ast}} \;\middle|\; Q \in f^{-1}(V(S))
\right\}.
\end{align*}
Since $f$ is continuous, $f^{-1}(V(S)) = V(K)$ for some subhypermodule
$K \leq M$. Hence
\[
(f^{\ast})^{-1}
\left(
V\!\left(\frac{S}{\theta^{\ast}}\right)
\right)
=
V\!\left(\frac{K}{\theta^{\ast}}\right),
\]
which is closed. Therefore $f^{\ast}$ is continuous.

\medskip

Finally, identity maps are preserved and for composable morphisms
$f,g$ we have
\[
(g \circ f)^{\ast} = g^{\ast} \circ f^{\ast}.
\]
Hence $F$ is a functor.
\end{proof}

\end{document}